\documentclass[9pt,letterpaper]{IEEEtran4PSCC/IEEEtran4PSCC}
\ifCLASSINFOpdf
\else
\fi
\usepackage{rotating}
\usepackage[cmex10]{amsmath}
\usepackage{psfrag}
\usepackage{multicol}
\usepackage{epstopdf}
\usepackage{graphicx,cite}
\usepackage{subfigure,textcomp,gensymb,url}
\usepackage[usenames, dvipsnames]{color}
\usepackage{xcolor}
\usepackage{algorithm}
\usepackage{algpseudocode}%
\usepackage[inline]{enumitem}

\newcommand{\comments}[1]{}
\usepackage{anysize, booktabs, setspace}
\usepackage[inline]{enumitem}
\marginsize{0.7in}{0.7in}{0.37in}{0.05in}
\usepackage{amsfonts, amssymb, amsthm, eucal}
\usepackage{cleveref, nomencl, blindtext, mathtools, etoolbox}
\renewcommand\nomgroup[1]{%
  \item[\bfseries
  \ifstrequal{#1}{A}{Sets and Parameters}{%
  \ifstrequal{#1}{B}{Parameters}{%
  \ifstrequal{#1}{C}{Binary variables for Diesel Generators}{%
    \ifstrequal{#1}{D}{Continuous variables}{%
  \ifstrequal{#1}{E}{Continuous variables for Diesel Generators}{%
  \ifstrequal{#1}{F}{Continuous variables for Battery Resources}{%
   \ifstrequal{#1}{G}{Other(Line / Bus) Variables}{
    \ifstrequal{#1}{H}{Scenario Variables}{}}}}}}}}]%
  }
\immediate\write18{makeindex \jobname.nlo -s nomencl.ist -o \jobname.nls}
\makenomenclature
\newcommand{\fr}[1]{\mathfrak{#1}}
\newcommand{\bs}[1]{\boldsymbol{#1}}
\newcommand{\wh}[1]{\widehat{#1}}

\begin{document}

\title{Tight Piecewise Convex Relaxations for\\ Global Optimization of Optimal Power Flow}

%% To specify the authors when (number of affiliations <= 2)
% \author{
% \IEEEauthorblockN{Author n.1 Name per Affiliation A\\ Author n.2 Name per Affiliation A}
% \IEEEauthorblockA{(Affiliation A) Department Name of Organization \\
% Name of the organization, acronyms acceptable\\
% City, Country\\
% \{email author n.1, email author n.2\}@domain (if desired)}
% \and
% \IEEEauthorblockN{Author n.1 Name per Affiliation B\\ Author n.2 Name per Affiliation B}
% \IEEEauthorblockA{(Affiliation B) Department Name of Organization \\
% Name of the organization, acronyms acceptable\\
% City, Country\\
% \{email author n.1, email author n.2\}@domain (if desired)}
% }

\author{
   Mowen Lu$^{\dag}$, Harsha Nagarajan$^\ddag$, Russell Bent$^*$, Sandra D. Eksioglu$^{\dag}$, Scott J. Mason$^{\dag}$ \\
$^\dag$ Department of Industrial Engineering, Clemson University, SC, United States.\\
$^\ddag$ Theoretical Division (T-5), Los Alamos National Laboratory, NM, United States.  Contact: harsha@lanl.gov \\
$^*$Analytics, Intelligence and Technology Division (A-1), Los Alamos National Laboratory, NM, United States.
}

%% To specify the authors when (number of affiliations > 2)
% \author{\IEEEauthorblockN{Author n.1\IEEEauthorrefmark{1},
% Author n.2\IEEEauthorrefmark{2},
% Author n.3\IEEEauthorrefmark{3}, 
% Author n.4\IEEEauthorrefmark{3} and
% Author n.5\IEEEauthorrefmark{4}}
% \IEEEauthorblockA{\IEEEauthorrefmark{1} Department Name of Organization A\\
% Name of the organization A,
% Address A\\ Emails if wanted}
% \IEEEauthorblockA{\IEEEauthorrefmark{2} Department Name of Organization B\\
% Name of the organization B,
% Address B\\ Emails if wanted}
% \IEEEauthorblockA{\IEEEauthorrefmark{3} Department Name of Organization C\\
% Name of the organization C,
% Address C\\ Emails if wanted}
% \IEEEauthorblockA{\IEEEauthorrefmark{4}Department Name of Organization D\\
% Name of the organization D,
% Address D\\ Emails if wanted}
% }

\maketitle

\begin{abstract}
Since the alternating current optimal power flow (ACOPF) problem was introduced in 1962, developing efficient solution algorithms for the problem has been an active field of research. In recent years, there has been increasing interest in convex relaxations-based solution approaches that are often tight in practice. Based on these approaches, we develop tight piecewise convex relaxations with convex-hull representations, an adaptive, multivariate partitioning algorithm with bound tightening that progressively improves these relaxations and, given sufficient time, converges to the globally optimal solution. We illustrate the strengths of our algorithm using benchmark ACOPF test cases from the literature. Computational results show that our novel algorithm reduces the best-known optimality gaps for some hard ACOPF cases.
\end{abstract}


%----------------------------------;
%    Include all sections here     ;
%----------------------------------;
%!TEX root = ../pscc2018.tex

%+++++++++++++++++++++++++++;
%  Section: Introduction      ;
%+++++++++++++++++++++++++++;
\section{Introduction}
\label{Sec:intro}
The optimal power flow problem is one of the most fundamental optimization problems for the economic and reliable operation of electric power systems. The ACOPF is a cost minimization problem with equality and inequality constraints setting bus voltage, line flows and generator dispatch. It was first introduced in 1962 and has been formulated in various forms over years, e.g., the polar power-voltage formulation and the rectangular power-voltage formulation \cite{cain2012history}. Since the introduction of the ACOPF problem \cite{carpentier1962contribution}, the development of efficient solution technique for the ACOPF has remained an active field of research. The main challenges associated with solving the ACOPF include: a) non-convex and nonlinear mathematical models of ac physics, b) large-scale power grids, and c) limited computation time available in real-time dispatch applications. A fast and robust algorithm for obtaining high-quality and lower-cost dispatch solutions could improve the operational performance of power systems and save billions of dollars per year \cite{cain2012history}.

Various solution algorithms, ranging from heuristic to convex relaxation, have been studied within literature \cite{mccormick1976computability,hijazi2017convex}. Without losing generality, the ACOPF literature can be roughly categorized into three main research directions: a) finding locally optimal solutions quickly, b) deriving strong, convex relaxations and c) proving global optimality. In the local search literature, local solvers based on primal-dual interior point methods or sequential linearization heuristics are used to find feasible solutions efficiently without any guarantees on solution quality. In the relaxation literature, recent work has focused on deriving convex relaxations that produce tight lower bounds. These include SDP-based methods \cite{lavaei2012zero}, Quadratic relaxations (QC) \cite{hijazi2017convex} and Second-order-cone (SOC) relaxations \cite{kocuk2017matrix}. The performances of these existing methods have been evaluated and tested on well established power system test cases (i.e., Matpower and NESTA Test Cases) \cite{coffrin2014nesta}. Though the relaxation approaches have empirically yielded strong lower bounds, there remain examples where the lower bounds are weak (e.g., case5, nesta\_case30\_fsr\_api, nesta\_case118\_ieee\_api, etc.) \cite{hijazi2017convex}. Finally, there have been recent efforts focused on obtaining globally optimal solutions through SOC and SDP-based relaxations \cite{kocuk2017matrix, gopalakrishnan2012global} using standard spatial branch-\&-bound (sBB) approaches. Here, we focus on closing the optimality gaps on remaining hard instances by improving the QC relaxations, improved branching strategies, and leveraging high quality locally optimal solutions. 
    
The focus of this paper is a novel approach for globally optimizing the ACOPF. Our work is built on an adaptive multivariate partitioning algorithm (AMP) proposed in \cite{nagarajan2016tightening,nagarajan2017adaptive}. The approach is based on a two-stage algorithm that uses sBB-like methods tailored to OPF problems. In the first stage, we apply sequential bound-tightening techniques to the voltage and phase-angle variables and obtain tightest possible bounds by solving a sequence of convex problems \cite{coffrin2015strengthening,chen2016bound}. The second stage adaptively partitions convex envelopes of the ACOPF into piecewise convex regions around best-known local feasible solutions. 
This approach exploits the observation that local solutions to standard benchmark instances are already very good. 
%In comparison, standard sBB methods in global optimization partition envelopes uniformly \cite{hasan2010piecewise}. 
%Though this approach has many appealing properties, in the worst case this approach creates many partitions in unproductive areas of the search space. In these situations there is a combinatorial explosion in binary variables used to control the choice of partitions and leads to very slow convergence. 
Our recent results on generic mixed-integer nonlinear programs suggest that refining variable domains adaptively around best-known feasible solutions can dramatically speed up the convergence to global optimum \cite{ nagarajan2017adaptive,wu2017adaptive}.

This paper makes three key contributions to solving the ACOPF problem. The first contribution develops an efficient partitioning scheme for tightening relaxations.  In multilinear relaxations, many approaches build uniform, piecewise relaxations via univariate or bivariate partitioning \cite{hasan2010piecewise}. One drawback of such approaches is that a large number of partitions may be needed to attain global optimum. Thus, these approaches are often restricted to small problems. To address inefficiencies of these approaches, we develop an adaptive tightening algorithm with non-uniform partitions, where we selectively partition convex envelopes that heuristically appear to tighten the relaxations.

%However, none of these approaches give attention to reducing the number of variables to partition. Partitioning all the variables in multi-linear terms for large-sized problems is time consuming. Instead, we only partition a subset of variables based on the gaps measured by their local feasible solutions with the initial relaxed solutions. This heuristic algorithm is inspired by the observation that local solutions to standard benchmark instances are already very good.

Our second contribution lies in applying well-known ideas for deriving the tightest possible convex relaxations (convex-hulls) for multilinear functions as they play a crucial role for developing efficient global optimization approaches. Multilinears (up to trilinear) appear in the polar form of the ACOPF, for which there has been a recent development in developing strong convex quadratic relaxations (QC) \cite{coffrin2016qc,hijazi2017convex}. However, these relaxations employ recursive McCormick envelopes to handle the trilinear terms, which rarely capture their convex  hulls. In the optimization literature, specifically for a trilinear function, \cite{meyer2004trilinear} (Meyer-Floudas envelopes) describes the convex hull by deriving all it's facets, and for a generic multilinear function, convex hull is typically formulated as a convex combination of the extreme points of the function \cite{Rikun1997}. Owing to the simplicity of the latter idea from the global optimization perspective, we further strengthen the QC relaxations to obtain tighter lower bounds for the OPF problem. Also, we develop tight piecewise convex relaxations of the convex-hull representation for trilinear and quadratic functions by extending the ideas of approximating univariate/bivariate functions  \cite{padberg2000approximating,vielma2015mixed}. To the best of our knowledge, this is the first paper which applies tight representations of piecewise convex relaxations for ACOPF. 

The third contribution of this paper is a novel algorithm that combines iterative partitioning (using piecewise relaxations) with bound tightening to globally solve the ACOPF (given sufficient time). This algorithm first applies ``optimality-based bound tightening" by solving a sequence of min. and max. problems on voltage and phase angle difference variables \cite{puranik2017domain, coffrin2015strengthening, chen2016bound}. 
%At this step, we reduce the number of variables being tightened by exploiting implicit relationships between variables.
Second, the algorithm iteratively partitions the convex envelopes to tighten the lower bound. Simultaneously, the algorithm updates locally feasible solutions based to tighten the upper bound. The combination of tightening the upper and lower bound yields an algorithm akin to sBB that determines the globally optimal solution.

This paper is organized as follows: In section \ref{Sec:form}, we first revisit the ACOPF problem and the improved QC relaxation with convex hull representation. Then, we introduce piecewise relaxations in section \ref{sec:piecewise}. Section \ref{sec:global} discusses our global optimization algorithm. In section \ref{Sec:results}, we evaluate the proposed improvements on various hard test cases and conclude the paper in section \ref{Sec:conc}.
%!TEX root = ../pscc2018.tex

%+++++++++++++++++++++++++++;
%  Section: Formulation      ;
%+++++++++++++++++++++++++++;
\section{AC Optimal Power Flow Problem}
\label{Sec:form} 

\noindent
\textbf{Nomenclature} \\
\textbf{Sets and Parameters} \\
{$\cal{N}$} - {set of nodes (buses)} \\
{$\cal{G}$} - {set of generators} \\
{${\cal{G}}_i$} - {set of generators at bus $i$} \\
{$\cal{E}$} - {set of \textit{from} edges (lines)} \\
{$\cal{E}^{R}$} - {set of \textit{to} edges (lines) } \\
{$\bs{c}_0, \bs{c}_1, \bs{c}_2$} - {generation cost coefficients} \\
{$\bs{i}$} - {imaginary number constant} \\
{$\bs{Y}_{ij} = \bs{g}_{ij} +\bs{i}\bs{b}_{ij}$} - {admittance on line $ij$}\\
{$\bs{S}_{i}^{\bs{d}} = \bs{p}_{i}^{\bs{d}} + \bs{i}\bs{q}_{i}^{\bs{d}}$} - {AC power demand at bus $i$} \\
{$\overline{\bs{S}}_{ij}$} - {apparent power limit on line $ij$} \\
{$\underline{\bs{\theta}}_{ij}, \overline{\bs{\theta}}_{ij}$} - {phase angle difference limits on line $ij$} \\
{$\bs{\theta}_{ij}^{\bs{M}}$} - {$\max(\mid \underline{\bs{\theta}}_{ij}\mid, \mid\overline{\bs{\theta}}_{ij}\mid)$ on line ij} \\
{$\underline{\bs{v}}_i, \overline{\bs{v}}_i$} - {voltage magnitude limit at bus $i$} \\
{$ \underline{\bs{S}}_{i}^{\bs{g}}$, $\overline{\bs{S}}_{i}^{\bs{g}}$} - {power generation limit at bus $i$} \\
{$\fr{R}(\cdot)$} - {real part of a complex number} \\
{$\fr{T}(\cdot)$} - {imaginary part of a complex number} \\
{$(\cdot)^*$} - {hermitian conjugate of a complex number} \\
{$\mid\cdot\mid, \angle\cdot$} - {magnitude, angle of a complex number} \\

\noindent
\textbf{Continuous variables} \\
{${V_i} = v_i\bs{e}^{\bs{i}\theta_i}$} - {AC voltage at bus $i$} \\
{$\theta_{ij} = \angle V_i - \angle V_j$} - {phase angle difference on line $ij$} \\
{$W_{ij}$} - {AC voltage product on line $ij$, i.e., $V_iV^*_j$} \\
{$S_{ij} = p_{ij} +\bs{i}q_{ij}$} - {AC power flow on line $ij$} \\
{$S^g_i= p^g_i+\bs{i}q^g_i$} - {AC power generation at bus $i$} \\
{$l_{ij}$} - {current magnitude squared on line $ij$} \\

In this paper, constants are typeset in bold face. In the AC power flow equations, the primitives, $V_i$, $S_{ij}$, $S^g_i$, $\bs{S}_{i}^{\bs{d}}$ and $\bs{Y}_{ij}$ are complex quantities. Given any two complex numbers (variables/constants) $z_1$ and $z_2$, $z_1 \geqslant z_2$ implies $\fr{R}({z}_1) \geqslant \fr{R}({z}_2)$ and $\fr{T}({z}_1) \geqslant \fr{T}({z}_2)$. $|\cdot|$ represents absolute value when applied on a real number. Statement $A \wedge B$ is true \textit{iff} $A$ and $B$ are both true; else it is false. $\langle f(\cdot)\rangle^R$ represents the constraints corresponding to the convex relaxation of function $f(\cdot)$.

This section describes the mathematical formulation of the ACOPF problem using the \textit{polar} formulation. Here, a power network is represented as a graph, $(\cal{N},\cal{E})$, where $\cal{N}$ and $\cal{E}$ are the buses and transmission lines, respectively. Generators are connected to buses where ${\cal{G}}_i$ are the generators at bus $i$. We assume that there is power demand (load) at every bus, some of which is zero. The optimal solution to the ACOPF problem minimizes generation costs for a specified demand and satisfies engineering constraints and power flow physics. More formally, the ACOPF problem is mathematically stated as:

\begin{subequations}
% \color{red}
\label{eq:AC_polar}
\allowdisplaybreaks
%\small
\begin{align}
&\label{objective}\bs{\cal{P}}:=\min \sum_{i\in \cal{G}} \bs{c}_{2i}(\fr{R}(S_i^g)^2) + \bs{c}_{1i}\fr{R}(S_i^g) + \bs{c}_{0i}\\
% &\text{s.t. } \nonumber\\
&\label{s_balance}\textit{s.t.} \ \  \sum_{k \in {\cal{G}}_i} S_k^g - \bs{S}_{i}^{\bs{d}} = \sum_{(i,j)\in\cal{E} \cup \cal{E^R}}S_{ij} \ \ \hspace{15pt}\forall i \in \cal{N} \\
&\qquad\label{sij}S_{ij}=\bs{Y}^{*}_{ij}W_{ii} - \bs{Y}^*_{ij}W_{ij} \hspace{15pt}\forall (i,j) \in \cal{E} \\
&\qquad\label{sji}S_{ji}=\bs{Y}^{*}_{ij}W_{jj} - \bs{Y}^{*}_{ij}W^*_{ij}\hspace{15pt}\forall (i,j) \in \cal{E}\\
&\qquad\label{wii} W_{ii} = |V_i|^2 \ \  \forall i \in \cal{N}\\
&\qquad\label{wij} W_{ij} = V_iV_j^* \ \  \forall (i,j) \in \cal{E}\\
& \qquad \label{theta} \underline{\bs{\theta}}_{ij} \leqslant \angle V_i - \angle V_j \leqslant \overline{\bs{\theta}}_{ij} \ \ \forall (i,j) \in \cal{E}\\
&\qquad\label{v_limit}\underline{\bs{v}}_i \leqslant |V_i| \leqslant \overline{\bs{v}}_i \ \  \forall i \in \cal{N} \\
&\qquad\label{g_cap} \underline{\bs{S}}_{i}^{\bs{g}} \leqslant S^g_i \leqslant \overline{\bs{S}}_{i}^{\bs{g}} \ \  \forall i \in \cal{G} \\
&\qquad\label{s_cap}|S_{ij}| \leqslant \overline{\bs{S}}_{ij} \ \  \forall (i,j) \in \cal{E}\cup \cal{E^R}
\end{align}
\end{subequations}

In formulation \eqref{eq:AC_polar}, the convex quadratic objective \eqref{objective} minimizes total generator dispatch cost. 
Constraint \eqref{s_balance} corresponds to the nodal power balance at each bus, e.g. Kirchoff's current law. Constraints \eqref{sij} through \eqref{wij} model the AC power flow on each line in complex number notation. Constraint \eqref{theta} limits the phase angle difference on each line. Constraint \eqref{v_limit} limits the voltage magnitude at each bus. Constraint \eqref{g_cap} restricts the apparent power output of each generator. Finally, constraint \eqref{s_cap} restricts the total electric power transmitted on each line. For simplicity, we omit the details of constant bus shunt injections, transformer taps, phase shifts, and line charging, though we include them in the computational studies.  

The ACOPF is a hard, nonconvex problem \cite{bienstock2015strong} where the source of \textit{nonconvexity} is in constraints \eqref{wii} and \eqref{wij}, which reduce to: 
\begin{subequations}
\begin{align}
W_{ii} &= v_i^2\\
\label{eq:wcs}\fr{R}(W_{ij}) &= v_iv_j\cos(\theta_{ij})  \\
\label{eq:wsn}\fr{T}(W_{ij}) &= v_iv_j\sin(\theta_{ij}) 
\end{align}
\label{eq:w_nonconv}
\end{subequations}
\vspace{-0.5cm}

To address the nonconvexities, we first summarize a state-of-the-art convex relaxation with some enhancements and then derive tighter piecewise relaxations. 

\section{Convex Quadratic Relaxation of the ACOPF}
\label{sec:convex}

In this section, we discuss the features of the convex quadratic (QC) relaxation of the ACOPF \cite{coffrin2015strengthening,hijazi2017convex}. Though there are numerous other relaxations in the literature, we adopt the QC relaxations as it has been observed to be empirically tight, computationally stable, and efficient. We further tighten the sinusoidal relaxations of \cite{coffrin2015strengthening,hijazi2017convex} for certain conditions and introduce the tightest possible convex relaxations for trilinear functions. 

\noindent
\underline{\textbf{Quadratic function relaxation}} \\
Given a voltage variable, $v_i \in [\underline{\bf v}_i,\overline{\bf v}_i]$, the tightest convex envelop is formulated with a lifted variable, $\wh{w}_i \in \langle v_i^2\rangle^{R}$, where 

\begin{subequations}
\label{eq:MC-q}
%\small
\allowdisplaybreaks
\begin{align}
&\wh{w}_i \geqslant v_i^2 \\
&\wh{w}_i \leqslant (\bs{\underline{v}}_{i}+\bs{\overline{v}}_{i})v_i - \bs{\underline{v}}_{i}\bs{\overline{v}}_{i}
% &\underline{v}_i \leqslant v_i \leqslant \overline{v}_i
\end{align}
\end{subequations}%

\noindent
\underline{\textbf{Cosine function relaxation}} \\
Under the assumptions that $|\underline{\bs{\theta}}_{ij}|$ is not always equal to $|\overline{\bs{\theta}}_{ij}|$ and  $\bs{\theta}_{ij}^{\bs{M}} \leqslant \pi/2$, the convex quadratic envelope of the cosine function is formulated with a lifted variable $\wh{cs}_{ij} \in \langle \cos(\theta_{ij})\rangle^{R}$ where 
\begin{subequations} \label{cosine}
\allowdisplaybreaks
\label{eq:cos}
%\footnotesize
\allowdisplaybreaks
\begin{align}
&\wh{cs}_{ij} \leqslant 1 - \frac{1- \cos(\bs{\theta}_{ij}^{\bs{M}})}{(\bs{\theta}_{ij}^{\bs{M}})^2}(\theta_{ij}^2) \\
&\wh{cs}_{ij} \geqslant \frac{\cos(\overline{\bs{\theta}}_{ij}) - \cos(\underline{\bs{\theta}}_{ij})}{\overline{\bs{\theta}}_{ij}-\underline{\bs{\theta}}_{ij}}(\theta_{ij} - \overline{\bs{\theta}}_{ij}) + \cos(\underline{\bs{\theta}}_{ij}) 
\end{align}
\end{subequations}%

\noindent
\underline{\textbf{Sine function relaxation}} \\ First, let the first-order Taylor's approximation of $\sin(\theta)$ at $\overline{\bs{\theta}}$ be
$$ f^{oa}(\overline{\bs{\theta}}) = \sin(\overline{\bs{\theta}}) + \cos(\overline{\bs{\theta}})(\theta-\overline{\bs{\theta}}).$$ 
Second, let the secant function between  $(\overline{\bs{\theta}},\sin(\overline{\bs{\theta}}))$ and $(\underline{\bs{\theta}},\sin(\underline{\bs{\theta}}))$ be
$$ f^{sec}(\underline{\bs{\theta}},\overline{\bs{\theta}}) = \frac{\sin(\overline{\bs{\theta}}) - \sin(\underline{\bs{\theta}})}{\overline{\bs{\theta}}-\underline{\bs{\theta}}}(\theta - \underline{\bs{\theta}}) + \sin(\underline{\bs{\theta}}).$$

\noindent
Finally, given bounds $ [\underline{\bs{\theta}},\overline{\bs{\theta}}]$, we use $\bs{\theta}^{\bs{\mu}} = \frac{(\bs{\underline{\theta}} + \bs{\overline{\theta}})}{2}$ to denote the midpoint of the bounds. 

Since the bounds on phase-angle differences can be non-symmetric, we derive convex relaxations of the sine function for three cases:\\
\noindent
\textbf{Case (a)}: When $(\underline{\bs{\theta}}_{ij} < 0) \wedge (\overline{\bs{\theta}}_{ij} > 0)$, the polyhedral relaxation, as described in \cite{hijazi2017convex}, is characterized with lifted variable $\wh{sn}_{ij} \in \langle \sin(\theta_{ij})\rangle^{R}$ where 
\begin{align} \label{case1}
%    \footnotesize
    \wh{sn}_{ij} \leqslant f^{oa}(\bs{\theta}_{ij}^{\bs{M}}/2), \ \ \wh{sn}_{ij} \geqslant f^{oa}(-\bs{\theta}_{ij}^{\bs{M}}/2) 
\end{align}

\noindent
\textbf{Case (b)}: When $(\underline{\bs{\theta}}_{ij} < 0) \wedge (\overline{\bs{\theta}}_{ij} < 0)$, we derive a tighter polyhedral relaxations that exploits a lack of an inflection point in the sine function. This relaxation is defined by: 
\begin{subequations} \label{case2}
%\footnotesize
\begin{align} 
    &\wh{sn}_{ij} \leqslant f^{sec}(\underline{\bs{\theta}}_{ij},\overline{\bs{\theta}}_{ij}), \\
    &\wh{sn}_{ij} \geqslant f^{oa}(\bs{\theta}_{ij}) \ \ \forall \bs{\theta}_{ij} \in \{\underline{\bs{\theta}}_{ij}, \bs{\theta}_{ij}^{\bs{\mu}}, \overline{\bs{\theta}}_{ij}\}.
\end{align}
\end{subequations}

\begin{figure}[htp]
   \centering
   \subfigure[Case (a)]{
   \includegraphics[scale=0.281]{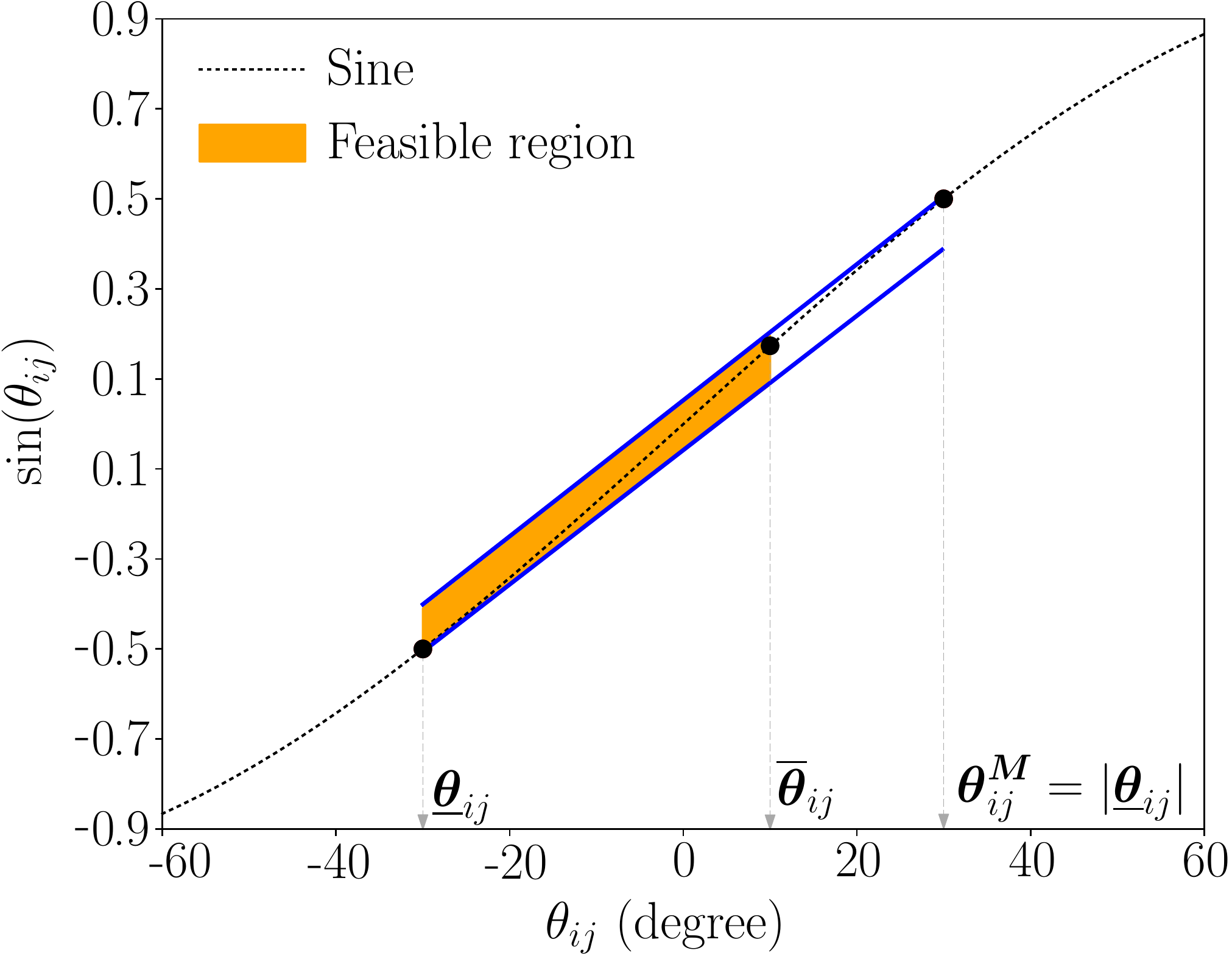}}       
   \subfigure[Case (c)]{
   \includegraphics[scale=0.28]{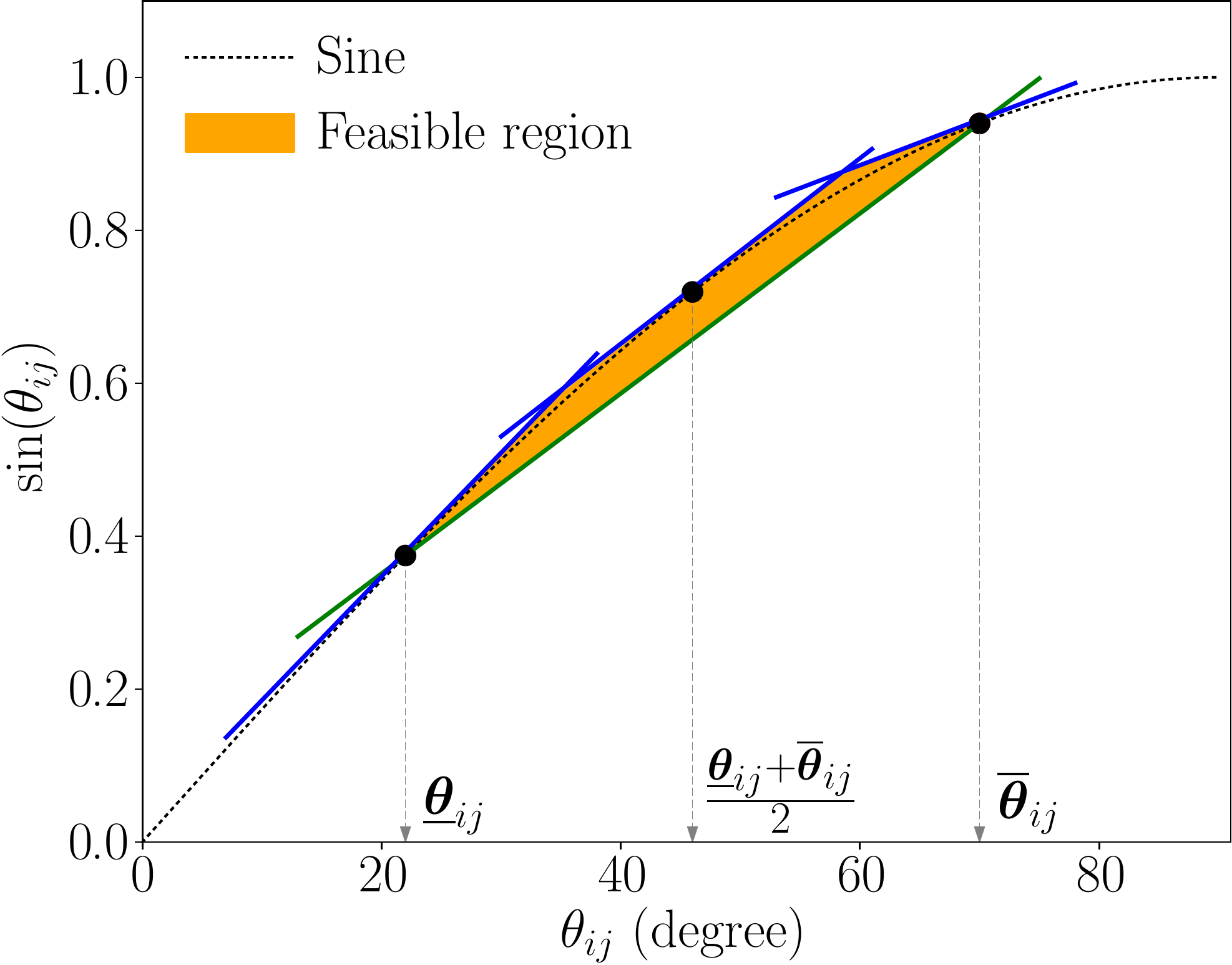}} 
   \caption{Polyhedral relaxations for sinusoidal function} 
   \label{Fig:sine_relax}
\end{figure}

\noindent
\textbf{Case (c)}: Like case (b), when $(\underline{\bs{\theta}}_{ij} > 0) \wedge (\overline{\bs{\theta}}_{ij} > 0)$, a tighter polyhedral relaxation for the sine function is: 
\begin{subequations} \label{case3}
%\footnotesize
\begin{align} 
    &\wh{sn}_{ij} \geqslant f^{sec}(\underline{\bs{\theta}}_{ij},\overline{\bs{\theta}}_{ij}), \\
    &\wh{sn}_{ij} \leqslant f^{oa}(\bs{\theta}_{ij}) \ \ \forall \bs{\theta}_{ij} \in \{\underline{\bs{\theta}}_{ij}, \bs{\theta}_{ij}^{\bs{\mu}}, \overline{\bs{\theta}}_{ij}\}.
\end{align}
\end{subequations}
A geometric visualization of the relaxations of the sine function for cases (a) and (c) are shown in Figure \ref{Fig:sine_relax}. 

\noindent
\underline{\textbf{Trilinear function relaxation}} \\
After introducing the lifted variables, $\wh{cs}_{ij}$ and $\wh{sn}_{ij}$ from above, the non-convex constraints in \eqref{eq:wcs} and \eqref{eq:wsn} become trilinear functions of the form $v_iv_j\wh{cs}_{ij}$ and $v_iv_j\wh{sn}_{ij}$. 
In the literature, bilinear McCormick relaxations are applied recursively to relax these trilinear functions \cite{mccormick1976computability,hijazi2017convex}, which rarely capture their convex hull. Instead, we relax the trilinear function based on the convex hull of the extreme points using techniques from \cite{Rikun1997}.

Given a trilinear function $\phi(x_1,x_2,x_3) = x_1x_2x_3$
with respective variable bounds $[\bs{\underline{x}_1}, \bs{\overline{x}_1}],[\bs{\underline{x}_2}, \bs{\overline{x}_2}],[\bs{\underline{x}_3}, \bs{\overline{x}_3}]$, the extreme points of $\phi(\cdot)$ are given by the Cartesian product
$ (\underline{\bs{x}}_1,\overline{\bs{x}}_1)\times (\underline{\bs{x}}_2,\overline{\bs{x}}_2) \times (\underline{\bs{x}}_3,\overline{\bs{x}}_3)=\langle \bs{\xi_1},\bs{\xi_2},\ldots, \bs{\xi_8} \rangle$. We use $\bs{\xi}^{i}_{\bs{k}}$ to denote the coordinate of $x_i$ in $\bs{\xi_k}$.
%Let $x_k^{EP},k=1,\ldots,2^3$ represent the elements of the set $\mathcal{X}^{EP}$.
The convex hull of the extreme points
%Given the extreme points, the \textit{tightest} convex relaxation (convex hull) 
of $x_1x_2x_3$ is then given by
%for the set 
%$$\small \mathcal{X} = \left\{(x_1,x_2,x_3,\wh{x}) \in ([\underline{\bs{x}}_i,\overline{\bs{x}}_i], \forall i=1,2,3) \times \bb{R} : \wh{x} = x_1x_2x_3 \right\},$$
%is given by 
\begin{subequations}\label{extreme_point}
%\footnotesize
\begin{align}
    &\sum_{k=1..8} \lambda_k = 1, \quad \lambda_k \geqslant 0, \quad \forall k=1,\ldots,8,\\
    &\wh{x} = \sum_{k=1..8} \lambda_k \phi(\bs{\xi_k}),\quad x_i = \sum_{k=1..8} \lambda_k \bs{\xi}^{i}_{\bs{k}}
\end{align}
\end{subequations}

\noindent The notation $\langle x_1,x_2,x_3\rangle^\lambda$ is used to denote the $\lambda$-based relaxation of a trilinear function as defined above. Thus, the relaxation of $x_1x_2x_3$ is stated as $\wh{x} = \langle x_1,x_2,x_3\rangle^\lambda$. We note that this formulation generalizes to any multilinear function and is equivalent to the standard McCormick relaxation for bilinear functions.  

\noindent
\underline{\textbf{Current-magnitude constraints}} \\
We also add the second-order conic constraints that connect apparent power flow on lines, ($S_{ij}$), with current magnitude squared variables, ($l_{ij}$) \cite{coffrin2016qc,hijazi2017convex}.
The complete convex quadratic formulation with \textit{tightest} trilinear relaxation of the ACOPF is then stated as:

% \vspace{-0.09cm}
%
\begin{subequations}
% \color{red}
\label{eq:QC}
\allowdisplaybreaks
%\small
\begin{align}
&\label{qc_objective}\bs{\mathcal{P}^{QC}}:=\min \sum_{i\in \cal{G}} \bs{c}_{2i}(\fr{R}(S_i^g)^2) + \bs{c}_{1i}\fr{R}(S_i^g) + \bs{c}_{0i}\\
% &\text{s.t. } \nonumber\\
& \textit{s.t.} \ \ \mathrm{Constraints} \ \eqref{s_balance}-\eqref{sji}, \\
& \qquad \mathrm{Constraints} \ \eqref{theta}-\eqref{s_cap}, \\
& \qquad W_{ii} = \wh{w}_i, \ \wh{w}_i \in \langle v_i^2\rangle^{R} \ \  \forall i \in \cal{N}\\
& \label{eq:wcs_wsn}\qquad \fr{R}(W_{ij}) = \wh{wcs}_{ij}, \ \fr{T}(W_{ij}) = \wh{wsn}_{ij},\ \  \forall (i,j) \in \cal{E} \\
% & \qquad\fr{T}(W_{ij}) = wsn_{ij}, \\
& \qquad \wh{wcs}_{ij} \in \langle v_iv_j\wh{cs}_{ij} \rangle^{\lambda}, \ \wh{wsn}_{ij} \in \langle v_iv_j\wh{sn}_{ij}\rangle^{\lambda}, \\
& \qquad \wh{cs}_{ij} \in \langle \cos(\theta_{ij})\rangle^{R}, \  \wh{sn}_{ij} \in \langle \sin(\theta_{ij})\rangle^{R}, \\
&\qquad\label{loss} S_{ij} + S_{ji} = Z_{ij}l_{ij} \ \ \forall (i,j) \in \cal{E} \\
&\qquad\label{SOC}|S_{ij}|^2 \leqslant W_{ii}l_{ij} \ \  \forall (i,j) \in \cal{E}.
\end{align}
\end{subequations}

% \input{Sections/Formulation_arx}
%!TEX root = ../pscc2018.tex

%+++++++++++++++++++++++++++;
%  Section: Algorithms      ;
%+++++++++++++++++++++++++++;

{\color{black}

\section{Piecewise Convex Relaxations}
\label{sec:piecewise}
One of the weaknesses of the convex quadratic relaxations described in section \ref{sec:convex} is that the relaxation is not tight when the bounds of the variables are wide.
To address this issue, recent work \cite{coffrin2015strengthening,chen2016bound} has developed approaches to tighten variable bounds, sometimes significantly.
However, there are still a few OPF instances 
with large optimality gaps. In this section, we focus on developing tighter piecewise convex relaxations for quadratic and trilinear functions. To the best of our knowledge, this is the \textit{first} paper that considers global optimization of ACOPF via piecewise polyhedral relaxations.    

\begin{figure}[htp]
   \centering
   \includegraphics[scale=0.8]{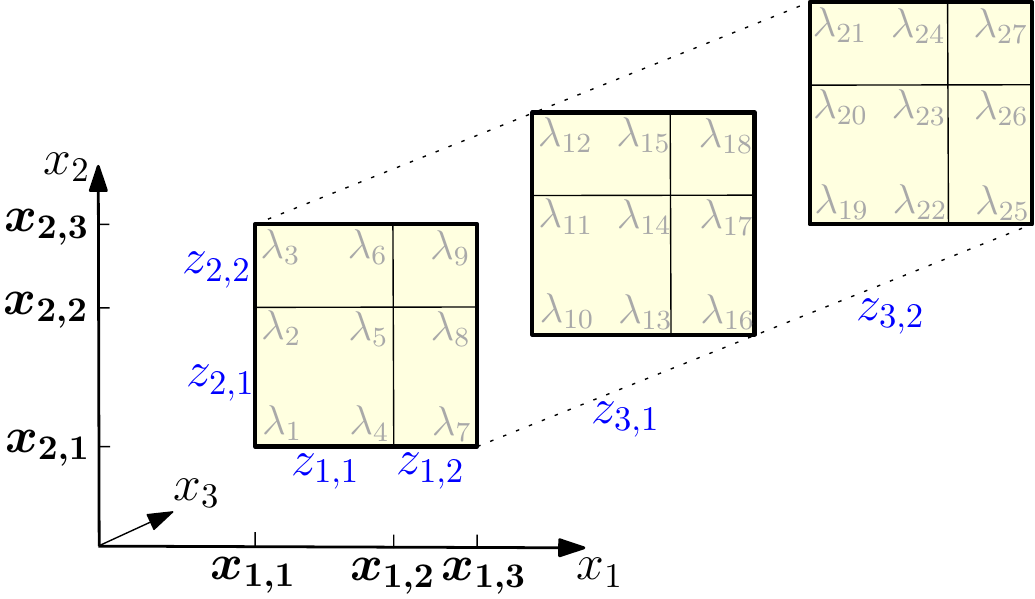}
   \caption{Partitioned variable domains for a trilinear function. Blue font indicates binary variables.}
   \label{fig:lambda}
\end{figure}

\noindent
\underline{\textbf{Piecewise trilinear functions}} \\ The basis for our piecewise relaxation is the piecewise linear approximation of continuous (mostly bivariate) separable functions that are solved with mixed integer programs (the $\lambda$-method) \cite{padberg2000approximating,lee2001polyhedral,vielma2015mixed}. 
Based on these approaches, we tighten the trilinear function relaxations  of $v_iv_j\wh{cs}_{ij}$ and $v_iv_j\wh{sn}_{ij}$ with a $\lambda$ formulation for modeling piecewise unions of polyhedrons. 

The piecewise relaxation partitions the variable domains by introducing discretization points within the variable bounds. Binary variables are used to control the selection of a partition and its associated local convex hull relaxation (in \eqref{extreme_point}) of the trilinear function.
%By assigning a binary variable per partition, it is used to control the partition that is active and thus the tighter convex hull relaxation (in \eqref{extreme_point}) associated with the active partition is applied to the trilinear function. 
As the number of partitions approaches $\infty$, the partitioning exactly represents the original trilinear function. Without loss of generality and for ease of exposition, we restrict the discussion of piecewise polyhedral relaxations for trilinear functions with \textit{two} partitions on every variable. 

\noindent
\textbf{Notation:} Let the discretization points of variable $x_i$ be $(\bs{\underline{x}_{i}} = \bs{x_{i,1}}) \leqslant \bs{x_{i,2}} \leqslant (\bs{x_{i,3}} = \bs{\overline{x}_{i}})$. Also, let $z_{i,1} \in \{0,1\}$ and $z_{i,2} \in \{0,1\}$ be the binary partition variables of $x_i$. We define $\mathcal{K}=\{1,..,27\}$ as the set of the first 27 positive integers. 
For $k\in \mathcal{K}$, $\bs{\xi}_{k}$ is used to denote the coordinates of the $k$th extreme point.
%Let $\bs{\xi}_{k}, k\in \mathcal{K}$ represent the set of co-ordinates of partition's extreme points given by $(\bs{x_{1,i}},\bs{x_{2,j}},\bs{x_{3,k}}),\forall i,j,k=1,2,3$. 
For every extreme point, $\bs{\xi}_{k}$, there is a nonnegative mulitiplier variable $\lambda_k$. The partitioned domains of the trilinear function are graphically illustrated in Fig. \ref{fig:lambda}. Using this notations, we now present SOS-II-type  constraints that model the piecewise union of polyhedrons of a trilinear function:
\begin{subequations}
%\footnotesize
\allowdisplaybreaks
\begin{align}
& \label{eq:sos2_a}\wh{x} = \sum_{k \in \mathcal{K}} \lambda_{k} \phi(\bs{\xi}_{k}), \quad x_i = \sum_{k \in \mathcal{K}} \lambda_{k} \bs{\xi}_{k}^{i},\\
& \label{eq:sos2_b}\sum_{k \in \mathcal{K}}\lambda_k=1, \quad \lambda_k \geqslant 0, \ \forall k \in \mathcal{K}, \\
& \label{eq:sos2_c}z_{i,1}, z_{i,2} \in \{0,1\}, \quad  z_{i,1} + z_{i,2}=1, \ \forall i=1,..,3, \\
& \label{eq:sos2_d}z_{1,1} \geqslant \sum_{k\in \widetilde{\mathcal{K}}} \lambda_k, \ \widetilde{\mathcal{K}}=\{1,2,3,10,11,12,19,20,21\}, \\
& \label{eq:sos2_e}z_{1,2} \geqslant \sum_{k\in \widetilde{\mathcal{K}}} \lambda_k, \ \widetilde{\mathcal{K}}=\{7,8,9,16,17,18,25,26,27\}, \\
& \label{eq:sos2_f}z_{2,1} \geqslant \sum_{k\in \widetilde{\mathcal{K}}} \lambda_k, \ \widetilde{\mathcal{K}}=\{1,4,7,10,13,16,9,22,25\}, \\
& \label{eq:sos2_g}z_{2,2} \geqslant \sum_{k\in \widetilde{\mathcal{K}}} \lambda_k, \ \widetilde{\mathcal{K}}=\{3,6,9,12,15,18,21,24,27\}, \\
& \label{eq:sos2_h}z_{1,1} + z_{1,2} \geqslant \sum_{k\in \widetilde{\mathcal{K}}} \lambda_k, \ \widetilde{\mathcal{K}}=\{4,5,6,13,14,15,22,23,24\},\\
& \label{eq:sos2_i}z_{2,1} + z_{2,2} \geqslant \sum_{k\in \widetilde{\mathcal{K}}} \lambda_k, \ \widetilde{\mathcal{K}}=\{2,5,8,11,14,17,20,23,26\},\\
& \label{eq:sos2_j}z_{3,1} \geqslant \sum_{k=1}^{9} \lambda_k, \quad z_{3,1} + z_{3,2} \geqslant \sum_{k=10}^{18} \lambda_k, \quad z_{3,2} \geqslant \sum_{k=19}^{27} \lambda_k.
\end{align}
\label{eq:sos2}
\end{subequations}

% \vspace{-0.2cm}
Constraints \eqref{eq:sos2_a}-\eqref{eq:sos2_c} model the convex combination of extreme points and disjunctions. They force one partition for each variable to be active. Constraints \eqref{eq:sos2_d}-\eqref{eq:sos2_j} enforce the adjacency conditions for the $\lambda$ variables. These resemble SOS-II constraints. The formulation in \eqref{eq:sos2} has many interesting polyhedral properties. For example, the projection of this polytope on to the space of $\{x_1,x_2,x_3,\wh{x}\}$ has integral extreme points. This is one of the primary reasons this formulation can be computationally efficient\footnote{The intersection of polytope \eqref{eq:sos2} with any other constraints may \textit{not} necessarily yield a new polytope with integral extreme points.} For further theoretical details, we delegate the reader to \cite{lee2001polyhedral,vielma2015mixed}.
% The convex hull of nonlinear terms are derived by a lambda method which is a well know method to model piecewise-linear functions \cite{vielma2010mixed, lee2001polyhedral}. We next formally define the convex hull of $x_1x_2x_3$, $conv(x_1,x_2,x_3)$ as projected (Proj) into the $\lambda$ space, i.e.

% $$
% \small
% conv(x_1,x_2,x_3)= \underset{x_1,x_2,x_3,\wh{x}}{\mathrm{Proj}} \begin{cases}
% x_i \in [\underline{\bs{x}}_i,\overline{\bs{x}}_i], \forall i=1,2,3\\
% \mathrm{Eqs.} \;\; \ref{extreme_point}
% \end{cases}
% $$

\noindent
\underline{\textbf{Piecewise quadratic functions}}
% \vspace{-0.3cm}
\begin{figure}[htp]
   \centering
   \includegraphics[scale=0.82]{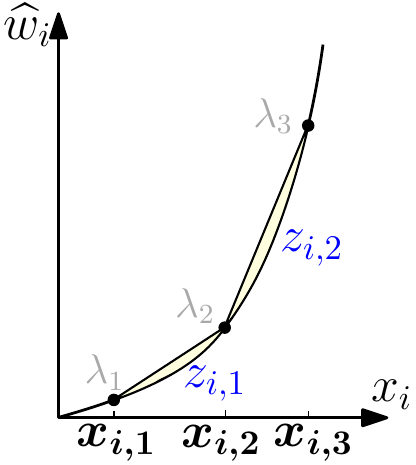}
   \caption{Piecewise quadratic regions. Blue font indicates binary variables.}
   \label{fig:quad_sos2}
\end{figure}
% \vspace{-0.15cm}

Using the same notation, the SOS-II constraints of piecewise quadratic constraints are modeled using a
piecewise union of convex quadratic regions (Fig. \ref{fig:quad_sos2}). 
To the best of our knowledge, this is the first time the piecewise quadratic regions of the voltage squared variables have been modeling using the $\lambda$ formulation.
\begin{subequations}
%\footnotesize
\begin{align}
& \label{eq:sos2q_a}\wh{w}_i \geqslant x_i^2,\quad \wh{w}_i \leqslant \sum_{k=1,2,3} \lambda_{k} \bs{x_i}^{2}_{k}, \quad x_i = \sum_{k=1,2,3} \lambda_{k} \bs{x_i}_{k},\\
& \label{eq:sos2q_b}\sum_{k=1,2,3}\lambda_k=1, \quad \lambda_k \geqslant 0, \ \forall k=1,2,3, \\
& \label{eq:sos2q_c}z_{i,1}, z_{i,2} \in \{0,1\}, \quad  z_{i,1} + z_{i,2}=1, \\
& \label{eq:sos2q_d} z_{i,1} \geqslant \lambda_1, \quad z_{i,2} \geqslant \lambda_3, \quad z_{i,1}+z_{i,2} \geqslant \lambda_2.
\end{align}
\label{eq:sos2_quad}
\end{subequations}
\noindent
\underline{\textbf{Strengthening valid inequalities}} \\
Though the formulations in \eqref{eq:sos2} and \eqref{eq:sos2_quad} are \textit{necessary and sufficient} to characterize the piecewise relaxations of trilinear and quadratic functions, we observed that the inclusion of the following (simple) valid constraints improved the computational performance of $\lambda$ formulations tremendously. For a given variable $x_i$ with two partitions (as described earlier), the constraint is as follows: 
\begin{align}
z_{i,1}\bs{x_{i,1}} + z_{i,2}\bs{x_{i,2}} \leqslant x_i \leqslant z_{i,1}\bs{x_{i,2}} + z_{i,2}\bs{x_{i,3}}
\label{eq:bnd}
\end{align}

Finally, for a given trilinear or quadratic function with a finite number of partitions, we use $\langle\cdot\rangle^{\lambda_p}$ to denote the piecewise $\lambda$-formulation of \eqref{eq:sos2} and \eqref{eq:sos2_quad}, respectively. The complete piecewise convex relaxation of the ACOPF is then stated as:
\begin{subequations}
% \color{red}
\label{eq:PQC}
\allowdisplaybreaks
%\small
\begin{align}
&\label{qc_objective}\bs{\mathcal{P}^{QC^{\lambda}}}:=\min \sum_{i\in \cal{G}} \bs{c}_{2i}(\fr{R}(S_i^g)^2) + \bs{c}_{1i}\fr{R}(S_i^g) + \bs{c}_{0i}\\
% &\text{s.t. } \nonumber\\
& \textit{s.t.} \ \ \mathrm{Constraints} \ \eqref{s_balance}-\eqref{sji}, \ \eqref{theta}-\eqref{s_cap}\\
& \qquad \mathrm{Constraints} \ \eqref{eq:wcs_wsn}, \ \eqref{loss}-\eqref{SOC}\\
& \qquad W_{ii} = \wh{w}_i, \ \wh{w}_i \in \langle v_i^2\rangle^{\lambda_p} \ \  \forall i \in \cal{N}\\
& \qquad \wh{wcs}_{ij} \in \langle v_iv_j\wh{cs}_{ij} \rangle^{\lambda_p}, \ \wh{wsn}_{ij} \in \langle v_iv_j\wh{sn}_{ij}\rangle^{\lambda_p}.
\end{align}
\end{subequations}
}

\section{Global optimization of ACOPF}
\label{sec:global}
\noindent
\underline{\textbf{Adaptive Multivariate Partitioning}} \\ To globally solve the ACOPF problem, we use the Adaptive Multivariate Partitioning (AMP) algorithm described in \cite{nagarajan2016tightening, nagarajan2017adaptive}. The key idea of AMP is that AMP leverages the observations that solutions based on relaxations to ACOPF are often tight in practice and that locally optimal solutions are also very good \cite{hijazi2017convex}. AMP iteratively introduces narrow partitions around these relaxations.

A high-level pseudo-code for AMP is given in Algorithm \ref{alg:AMP}. In this algorithm, we use the notation $\sigma$ to denote a solution to the ACOPF, $f(\sigma)$ to denote the objective value of $\sigma$, and $\sigma(x)$ to denote the assignment of variable $x$ in $\sigma$.
In Lines \ref{ac}-\ref{qc}, a feasible solution, $\overline\sigma$, and a lower bound, $\underline\sigma$, are computed. Here, the lower bound is computed without partitioning any variables.
In line \ref{bt}, we sequentially tighten the bounds of the voltage magnitude, ${v_i}$, and phase angle differences, $\theta_{ij}$, using optimization-based bound tightening (BT) \cite{coffrin2015strengthening,chen2016bound}. The new bounds on $\theta_{ij}$ are used to tighten $\widehat{cs}_{ij}$ and $\widehat{sn}_{ij}$ using the cases defined in equations \ref{cosine}-\ref{case3}. We would like to note that, though the importance of BT has been already observed for ACOPF problems, the bounds we obtain in this paper are tighter than in \cite{coffrin2015strengthening,chen2016bound} since formulation \eqref{eq:QC} is based on the convex-hull representation of trilinear functions.

\begin{figure}[htp]
   \centering
   \includegraphics[scale=0.35]{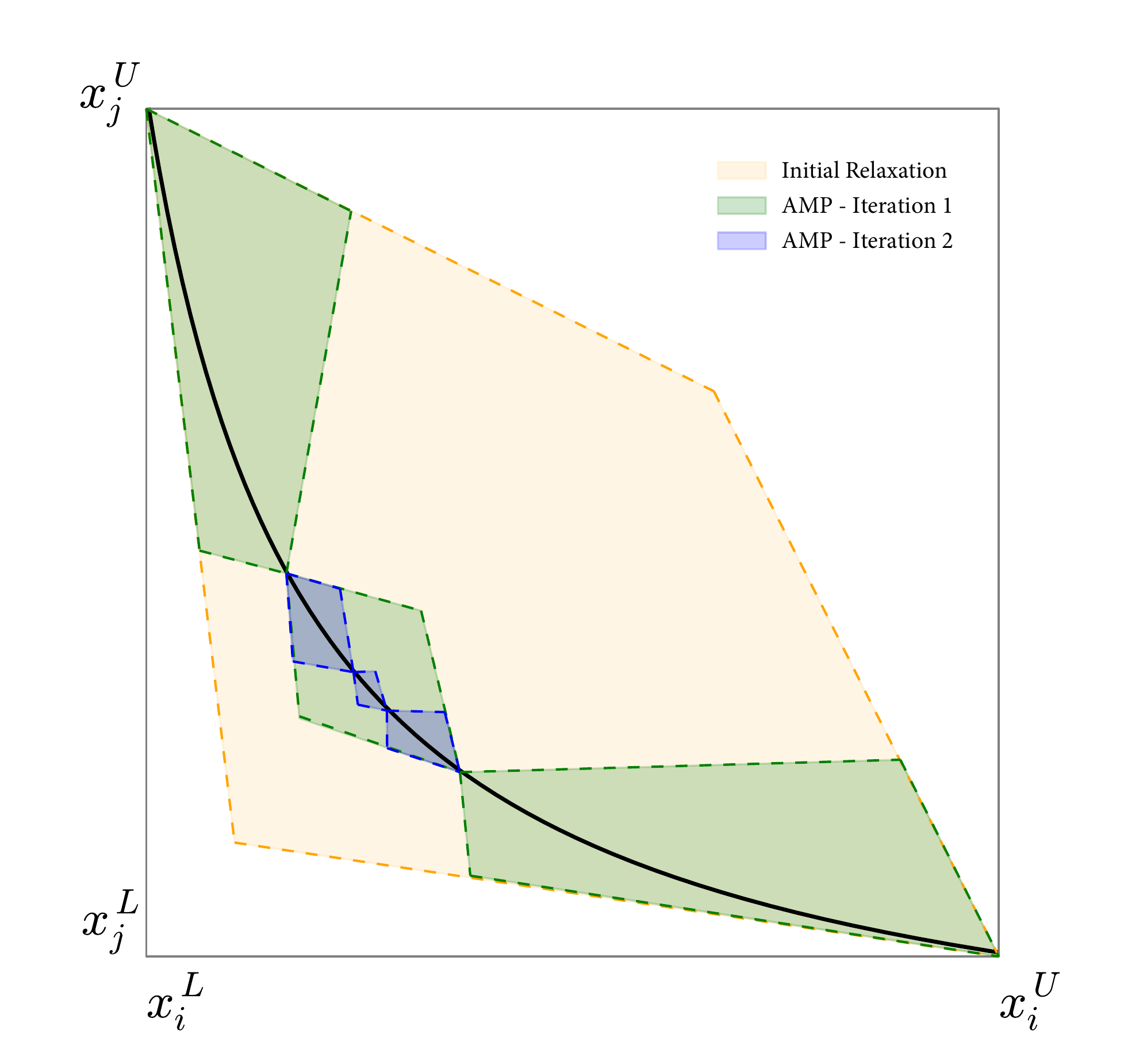}
   \caption{Example of variable partitioning in AMP for a bilinear function. Given an initial relaxation (tan region) and an initial feasible point (middle of the black curve), the function \textsc{InitializePartitions} creates a narrow partition around the feasible point (middle green region) and two wide partitions (outer green regions) around it. The size of the narrow partition is controlled by a user parameter, $\Delta$. The relaxation is iteratively tightened by the function \textsc{TightenPartitions} around relaxed solutions. This figure shows the partitioning of one iteration (blue regions).}
   \label{fig:region}
\end{figure}

Line \ref{heuristic_disc} describes our variable selection strategy for partitioning (discussed later). This is the main point of departure (and contribution) from the AMP algorithm discussed in \cite{nagarajan2017adaptive}.
Line \ref{partitionset} initializes the piecewise partitions of variables in $\mathcal{X}$ around those variables' assignments in $\overline\sigma$ (see Figure \ref{fig:region}). 
Line \ref{amp_1_iteration} then updates the lower bound using the piecewise relaxation.
Line \ref{amp:tighten} creates new partitions $\mathcal{X}$ around the updated lower bound (see Figure \ref{fig:region}). Lines \ref{amp:newlb}-\ref{amp:newub} then update the upper and lower bounds of the ACOPF. The upper bound (feasible solution) is updated by fixing the the ACOPF to the partition selected in $\underline\sigma$. and attempting to find a (better) feasible solution using a local solver. The process is repeated until the objective values of the upper and lower bound converge or a time-out criteria is met (Line \ref{beginwhile_amp}). The full details of all these procedures are discussed in \cite{nagarajan2017adaptive}. We note that this is the first time this algorithm has been applied to ACOPF represented in the polar form.

\begin{algorithm}[htp!]
\color{black}
% \footnotesize
%\scriptsize
\caption{Adaptive Multivariate Partitioning (AMP)}
\label{alg:AMP}
\begin{algorithmic}[1]
    \Function{AMP}{}

        \State\label{ac} $\overline{{\sigma}}  \gets$ \Call{Solve}{${\mathcal{P}}$}

        \State\label{qc} $\underline{{\sigma}}  \gets$ \Call{Solve}{${\mathcal{P}^{QC}}$}

        \State\label{bt} $\underline{\theta}_{ij}, \overline{\theta}_{ij}, \underline{v}_i, \overline{v}_i, 
        \underline{\widehat{cs}}_{ij}, 
        \overline{\widehat{cs}}_{ij}, 
        \underline{\widehat{sn}}_{ij},
        \overline{\widehat{sn}}_{ij}
        \gets $\Call{TightBounds}{${\overline{{\sigma}} }$}

       % This is the point of "newness" where % we decided which variables to partition
       \State\label{heuristic_disc} $\mathcal{X} \gets$\Call{SelectPartitionVariables}{$\overline{{\sigma}}, \underline{{\sigma}}$}

       \State\label{partitionset} ${\mathcal{I}} \gets $\Call{InitializePartitions}{$ \mathcal{X}, \mathcal{P}, \overline{{\sigma}}$}

       \State\label{amp_1_iteration} $\underline{{\sigma}} \gets$ \Call{Solve}{${\mathcal{P}^{QC^{\lambda}}}$}  
       
       \While{$\left(\frac{f(\overline{\sigma})-f(\underline{\sigma})}{f^(\underline{\sigma})}\geq \epsilon\right)$ and (Time $\leq$ Timeout)}\label{beginwhile_amp}

        \State \label{amp:tighten} ${\mathcal{P}^{QC^{\lambda}}} \gets $\Call{TightenPartitions}{$ x^{{\mathcal{I}}}, \mathcal{P}, \underline{\sigma}$}
        
       %\State $\underline{\sigma}' \gets \underline{\sigma}$

       \State \label{amp:newlb}$\underline{\sigma} \gets$\Call{Solve}{${\mathcal{P}^{QC^{\lambda}}}$}

       \State $\wh{\sigma} \gets$\Call{Solve}{$\mathcal{P}, \underline{\sigma}$}
       
       \If{$f(\wh{\sigma}) < f(\overline\sigma)$  }
       \State \label{amp:newub} $\overline\sigma \gets \wh{\sigma}$
       \EndIf

       \EndWhile\label{endwhile_amp}
       
        \State \Return $\underline{\sigma}, \overline{\sigma}$
    \EndFunction
\end{algorithmic}	
\end{algorithm}%

\noindent 
\underline{\textbf{Heuristic partition-variable selection}} \\
Algorithm \ref{alg:partitionedvariables} describes in more detail our variable selection strategy for partitioning. The original implementation of AMP selects a sufficient number of variables to ensure that all relaxed functions are partitioned and tightened. This approach ensures convergence to global optimality, but in practice, such an approach is computationally difficult to solve when the number of binary variables is large.
Here, we introduce a heuristic that limits the number of variables that are partitioned. While convergence is no longer assured, this heuristic can have a considerable impact on solution quality. 

Let $\mathcal{V} =\{v_i,\wh{cs}_{ij},\wh{sn}_{ij}, \ \forall i \in \mathcal{N}, ij \in \mathcal{E}\}$ represent the set of all variables appearing in nonconvex functions. In Algorithm \ref{alg:partitionedvariables}, Line \ref{delta} computes the difference between a variable's ($x_i$) assignment in an upper and a lower bound solution. Line \ref{sort} sorts the variables used in non-convex functions by increasing value of this difference. Line \ref{returnp} returns the the first $\alpha|\mathcal{S}|$ variables, where $\alpha$ is a user defined parameter between 0 and 1. This heuristic relies on an expectation that variables whose assignments from the non-partitioned convex relaxation are very different from the local feasible solution are indeed the variables that require further refinement in the relaxed space. Note that when $\alpha=1$ this heuristic reverts to the original AMP algorithm. We acknowledge that there can be numerous other branching strategies which can lead to better convergence of AMP, which we leave for the future work.

\begin{algorithm}[htp!]
\color{black}
% \small
%\footnotesize
\caption{Heuristic Partition-Variable Selection}
\label{alg:partitionedvariables}
\begin{algorithmic}[1]
    \Function{SelectPartitionVariables}{$\overline{{\sigma}},\underline{{\sigma}}$}

    %\For{$ x \in \{\mathbf{v}, \  \widehat{\mathbf{cs}}, \  \widehat{\mathbf{sn}}\} $} \label{forstart}
    
    \State $x_i^\delta \gets |\overline{\sigma}(x_i) - \underline{\sigma}(x_i)|,
    \ \ \forall i=1,\dots,|\mathcal{V}|$ \label{delta}
    
    %\EndFor\label{forend}
    
%    \State $\boldsymbol{\delta} \gets \{\Delta \sigma_x\}$
    
    \State $\mathcal{S} \gets \{ x_1, x_2,\dots,x_{|\mathcal{V}|} \} : x_i \in \mathcal{V}$ and $x_i^\delta \ge x_{i+1}^\delta$ \label{sort}
    %$|\overline{\sigma}(x_i) - \underline{\sigma}(x_i)| \ge |\overline{\sigma}(x_{i+1}) - \underline{\sigma}(x_{i+1})|$
    
     \State \Return $\bigcup x_i : i \le \alpha |\mathcal{S}|$ \label{returnp}
     
    \EndFunction
\end{algorithmic}	
\end{algorithm}%
%\vspace{-0.2cm}

% \input{Sections/Algo_arx}
%!TEX root = ../pscc2018.tex

%+++++++++++++++++++++++++++;
%  Section: Results      ;
%+++++++++++++++++++++++++++;
\section{Numerical Results}
\label{Sec:results}

In remainder of this paper, BT refers to bound tightening using the tightened convex quadratic relaxations introduced in section \ref{sec:convex}. QC$^{\mathrm{conv}}$ and BT-QC$^{\mathrm{conv}}$ correspond to improved QC relaxation and bound-tightening applied with convex-hull trilinear formulation \eqref{eq:QC}. BT-Lambda and BT-Lambda-$\alpha$ refer to Algorithm \ref{alg:AMP} with and without heuristic partitions (i.e., $\alpha = 100\%$ and $\alpha < 100\%$). The performance of these algorithms is evaluated on ACOPF test cases from NESTA 0.7.0 \cite{coffrin2014nesta} and are summarized in Table IV. 
These test cases were selected 
because the basic QC optimality gaps were larger than $1\%$ and hence hard for global optimization. Here, Ipopt 3.12.1 is used to find the feasible solution $\overline{\boldsymbol{\sigma}}$ in AMP.
The relaxed problems are solved 
using Gurobi 7.0.2 with default options and presolver switched on.

In Algorithm \ref{alg:AMP}, the value of $\epsilon$ and the ``time-out'' parameter were set to $0.01$ and $18000.0$ seconds, respectively. A 10800-second time-limit was imposed on the BT procedure. All implementations were done using Julia/JuMP  \cite{dunning2017jump}. Computational experiments were performed using the HPC Palmetto cluster at Clemson University with Intel Xeon E5-2670v2, 20 cores and 120GB of memory.

\subsection{Performance Comparison of Algorithms}

In Table IV, we summarize the computational performance of all the algorithms. The first two columns present the initial local feasible solution $\overline{\boldsymbol{\sigma}}$ and the initial lower bound $\underline{\boldsymbol{\sigma}}$, respectively. Columns three and four present the optimality gaps of \cite{coffrin2016qc} and \cite{kocuk2017matrix}, respectively. 
The remaining columns show the performance of BT, the piecewise relaxation, and heuristic variable partitioning.

Theses results show that the optimality gaps for BT-Lambda and BT-Lambda-$\alpha(60\%)$  are consistently smaller than  QC$^\mathrm{conv}$,
SDP gaps, and the best gaps in \cite{kocuk2017matrix}. In 24 out of 35 instances, the globally optimal solution is found (i.e., Gap $< 0.01\%$). Four of the sub-optimal solutions have gaps less than $0.2\%$. In addition, we observed that BT-Lambda-$\alpha(60\%)$ outperforms BT-Lambda on 29 out of 35 instances. For example,
on case73\_ieee\_rts\_sad, after
partitioning 60\% of the variables in nonlinear functions, the global optimal is found in 1.23 s ($58$x faster than all-variable partitioning).
The heuristic correctly identifies the subset of variables (fewer binary partitioning variables) required to prove optimality. 

The performance of the heuristic partitioning algorithm heavily depends on the choice of $\alpha$. Figure 5
%\ref{fig:alpha_sens} 
illustrates the computational gaps for different values of $\alpha$  on four test problems. In all four instances, ``GOpt" is found for all values of $\alpha$, but the changes in run times are dramatic.

\begin{figure}[htp!]
\centering
    \includegraphics[scale=0.31]{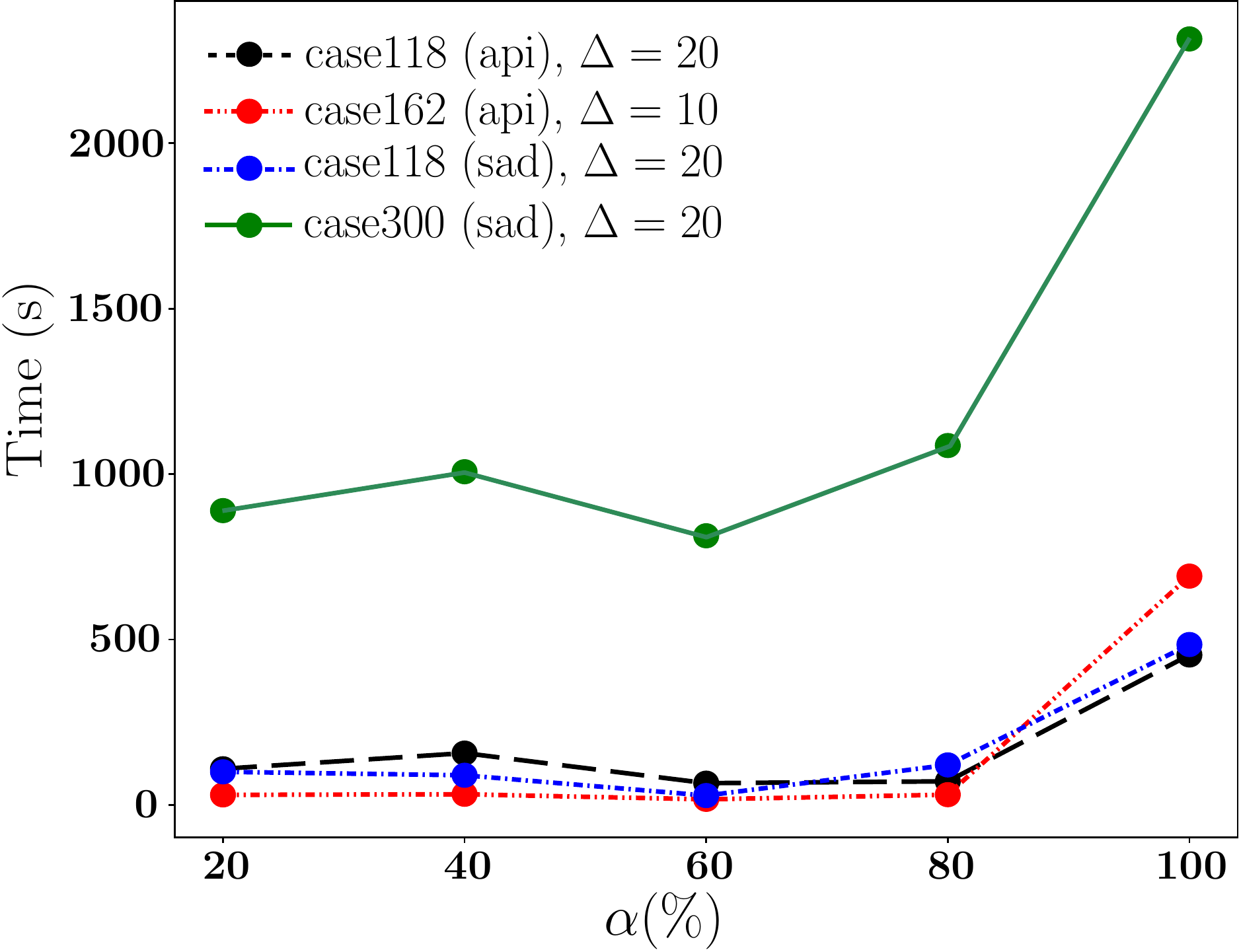}
    \label{fig:alpha_sens}
    \caption{$\alpha$ sensitivity analysis for run times.}
\end{figure}%

% \vspace{-0.2cm}
\noindent
\underline{\textbf{Improved QC relaxation gaps}} \\
\vspace{-0.5cm}
\begin{table}[htp]
    \centering
    \small
    \caption{QC relaxation gaps with trilinear functions relaxed using recursive McCormick (QC$^{\mathrm{rmc}}$) vs. convex-hull representation (QC$^{\mathrm{conv}}$).}
    \begin{tabular}{lrr}
    \toprule
         Instances & QC$^{\mathrm{rmc}}$ (\%) & QC$^{\mathrm{conv}}$ (\%) \\
         \midrule
         case3\_lmbd & 1.21 & 0.96 \\
         case30\_ieee & 15.64 & 15.20 \\
         \midrule
         case3\_lmbd\_api & 1.79 & 1.59 \\
         case24\_ieee\_rts\_api & 11.88 & 8.78 \\
         case73\_ieee\_rts\_api & 10.97 & 9.64 \\
         \midrule
         case3\_lmbd\_sad & 1.42 & 1.37 \\
         case4\_gs\_sad & 1.53 & 0.96 \\
         case5\_pjm\_sad & 0.99 & 0.77 \\
         case24\_ieee\_rts\_sad & 2.93 & 2.77 \\
         case73\_ieee\_rts\_sad & 2.53 & 2.38 \\
         case118\_ieee\_sad & 4.61 & 4.14 \\
\bottomrule
    \end{tabular}
    \label{tab:qc_gaps}
\end{table}%

Due to convex-hull relaxation of trilinear functions, the QC$^{\mathrm{conv}}$ gaps (highlighted in bold font in Table IV) outperformed the recursive McCormick approach used in \cite{coffrin2016qc}. Table \ref{tab:qc_gaps} tabulates a few instances in which the improvements in QC gaps were significant based on the proposed approach.  Even during bound-tightening, the best QC gaps using recursive McCormick (BT-QC$^{\mathrm{rmc}}$) for trilinear functions were 2.7\% and 11.8\% for case30\_fsr\_api and case118\_ieee\_api, respectively. However, applying the convex-hull representation (BT-QC$^{\mathrm{conv}}$) for the same instances significantly reduced the gaps to 0.35\% and 8.54\%, respectively.

\subsection{Analyses of ``Hard" Instances}
In Table IV,
%\ref{summary_tabel}, 
there are 7 instances where the optimality gaps of BT-Lambda and BT-Lambda-$\alpha(60\%)$ are larger than 1\%. These problems are hard because the tightened bounds remain weak.
Table \ref{tab:flip} shows that there are a large number of ${\theta}_{ij}$ variables whose bounds allow positive and negative values (flow could be in ether direction).
This is generally not the case for the other instances. We also noticed that BT-Lambda has better gaps than BT-Lambda-$\alpha(60\%)$ in 5/35 instances.
In this case, the heuristic incorrectly identifies variables whose partitioning is required to prove global optimality.
To illustrate this, Table \ref{tab:unparti_terms} counts the number of times a relaxed term does not have a variable that is partitioned (here, we use the lifted variable to refer to the term that is relaxed).

\noindent
\underline{\textbf{NCO instances}} AMP algorithm performed very well on 3/6 ``nco" instances, hard for nonconvex optimization, where AMP finds near global optimum solutions within the time limits. 

%
% \vspace{-0.6cm}
\begin{table}[htp]
    \centering
    \small
    \caption{\% of edges whose bounds are $(\bs{\underline{\theta}} < 0) \wedge (\bs{\overline{\theta}} > 0)$ after BT.}
    \begin{tabular}{lr}
    \toprule
         Instances & \% \\
         \midrule
        %  case5\_pjm & 33.3\\
         case240\_wecc & 77.6 \\
         case89\_pegase\_api & 51.5 \\
         case118\_ieee\_api & 32.4 \\
         case189\_edin\_api & 33.0 \\
         case240\_wecc\_sad & 68.1 \\
         case9\_bgm\_nco & 44.4 \\
         case39\_1\_bgm\_nco & 45.7\\
\bottomrule
    \end{tabular}
    \label{tab:flip}
\end{table}%

 \vspace{-0.6cm}
\begin{table}[htp!]
    \centering
    \small
    \caption{Relaxed terms without partitions}
    \begin{tabular}{lrrrrr}
\toprule
         Instances & $\widehat{w}$ & $\widehat{cs}$ & $\widehat{sn}$ & $\widehat{wcs}$ & $\widehat{wsn}$  \\
         \midrule
         case5\_pjm& 0 & 6 & 1 & 6 & 1 \\
         case240\_wecc\_sad & 22 & 258 & 94 & 264 & 124 \\
         case9\_na\_cao\_nco & 1 & 9 & 1 & 9 & 4 \\
         case9\_nb\_cao\_nco & 1 & 9 & 1 & 9 & 3 \\
         case14\_s\_cao\_nco & 0 & 19 & 3 & 19 & 3 \\
\bottomrule
    \end{tabular}
    \label{tab:unparti_terms}
\end{table}%

\begin{table*}[htp]
    \color{black}
    \centering
    \scriptsize
    \label{tab:summary_table}
    \caption{The summary of performance NESTA 0.7.0 ACOPF instances. Values under ``Gap'' and ``T,$ T^\alpha$'' are in \% and seconds, respectively. ``GOpt" refers to the global optimum when Gap $< 0.01\%$. ``TO" indicates time-out. ``--" indicates no solution is provided under ``SDP" and ``Best of \cite{kocuk2017matrix}" columns. ``--" under ``BT-Lambda" and ``BT-Lambda-$\alpha$" columns indicates that the BT-QC$^{\mathrm{conv}}$ already converged to global optimum. Here, $\Delta^*$ represents the best one from \{4, 6, 8, 10, 16\} and $\alpha = 60\%$. }
    \begin{tabular}{lrrrrrrrrrr}
\toprule
& & & & & 
& &\multicolumn{2}{c}{BT-Lambda} &
\multicolumn{2}{c}{BT-Lambda--$\alpha (60\%)$} \\
\cmidrule(lr){8-9}
\cmidrule(lr){10-11}
Instances & AC $(\$)$ & QC$^{\mathrm{conv}}$ (\%) & SDP (\%) & Best of  \cite{kocuk2017matrix} (\%) & BT$^{\mathrm{conv}}$ (s) & BT-QC$^{\mathrm{conv}}$ (\%) & ($\Delta^*$) Gap & $T$ & ($\Delta^*$) Gap & $T^{\alpha}$ \\
\midrule
case3\_lmbd & 5812.64 & \textbf{0.97} & 0.39 & 0.09 & 2.1 & GOpt & (16) GOpt & -- & (16) GOpt &--\\
case5\_pjm & 17551.90 & 14.54 & 5.22 & 0.10 & 0.9 & 6.73 & (6) GOpt & 303.7 & (6) \ 0.02 & 2023.6 \\
case30\_ieee& 204.97 & \textbf{15.20} & GOpt & -- & 48.1 & GOpt & (16) GOpt & -- & (16) GOpt & --\\
case118\_ieee & 3718.64 & \textbf{1.53} & 0.06 & 0.09 & 2235.0  & GOpt & (10) GOpt & -- & (20) GOpt & -- \\
case162\_ieee\_dtc & 4230.23 & 3.95 & 1.33 & 1.08 & 6292.3 & GOpt & (10) GOpt & -- & (10) GOpt & --\\
case240\_wecc & 75136.10 & \textbf{5.23} & -- & -- & 10800.0 & 2.90 &(6) \ 2.64  & TO  & (16) \ 2.56  & TO  \\
case300\_ieee & 16891.28 & 1.17 & 0.08 & 0.19 & 10800.0 & 0.01 & (20) GOpt & 2164.2 & (20) GOpt & 1091.8 \\
\midrule
case3\_lmbd\_api& 367.44 & \textbf{1.59} & 1.26 & 0.02 & 0.3 & GOpt & (16) GOpt & -- & (16) GOpt & -- \\
case14\_ieee\_api & 325.13 & 1.26 & GOpt & 0.09 & 11.7 & 0.06 & (8) GOpt & 0.4 & (16) GOpt &0.04\\
case24\_ieee\_rts\_api & 6426.65 & \textbf{8.79} & 1.45 & -- & 72.8 & GOpt & (8) GOpt & -- & (20) GOpt & --\\
case30\_as\_api & 570.08 & 4.63 & 0.00 & 0.06 & 47.1 & GOpt & (16) GOpt & -- & (16) GOpt & --\\
case30\_fsr\_api & 366.57 & 45.20 & 11.06 & 0.35 & 101.8 & 0.35 & (6) \ 0.12 & 75.5 & (8) \ 0.12 & 241.0\\
case39\_epri\_api& 7460.37 & 2.97 & GOpt & -- & 96.7 & GOpt & (16) GOpt & -- & (16) GOpt & --\\
case73\_ieee\_rts\_api& 19995.00 & \textbf{9.64} & 4.29 & -- & 1528.2 & 0.01 & (10) GOpt & 61.2 & (10) GOpt & 27.8\\
case89\_pegase\_api & 4255.44 & 19.83 & 18.11 & -- & 3997.3 & 16.21 & (6) 13.46  & TO  & (6) 13.08  & TO \\
case118\_ieee\_api& 10269.82 & 43.45 & 31.50 & 6.17 & 1515.6 & 8.54 & (6) \ 4.19  & TO & (6) \ 4.05 & TO \\
case162\_ieee\_dtc\_api &6106.86 & 1.25 & 0.85 & 1.03 & 6392.6 & GOpt & (6) GOpt & -- & (10) GOpt & -- \\
case189\_edin\_api & 1914.15 & 1.69 & 0.05 & 0.12 & 5110.9 & 0.04 & (8) \ 0.03 & 2351.2 & (8) \ 0.03 & 883.2  \\
\midrule
case3\_lmbd\_sad & 5959.33 & \textbf{1.38} & 2.06 & 0.03 & 0.2 & GOpt & (16) GOpt & -- & (20) GOpt & --\\
case4\_gs\_sad & 315.84 & \textbf{0.96} & 0.05 & -- & 0.3 & GOpt & (10) GOpt & -- & (10) GOpt & --\\
case24\_ieee\_rts\_sad  &76943.24 & \textbf{2.77} & 6.05 & -- & 44.7 & GOpt & (6) GOpt &-- & (16) GOpt & --\\
case29\_edin\_sad &41258.45 & 16.38 & 28.44 & 0.67 & 222.4 & GOpt & (16) GOpt &-- & (16) GOpt &--\\
case30\_as\_sad &897.49 & 2.32 & 0.47 & 0.08 & 53.8 & GOpt & (16) GOpt & -- & (16) GOpt &--\\
case30\_ieee\_sad &204.97 & \textbf{4.01} & GOpt & 0.08 & 50.4 & GOpt & (16) GOpt & -- & (16) GOpt &--\\
case73\_ieee\_rts\_sad  &227745.73 & \textbf{2.38} & 4.10 & -- & 1655.8 & GOpt & (6) GOpt & -- & (20) GOpt & --\\
case118\_ieee\_sad &4106.72 & \textbf{4.15} & 7.57 & 2.43 & 2211.4 & GOpt & (8) GOpt & -- & (20) GOpt & -- \\
case162\_ieee\_dtc\_sad &4253.51 & 4.27 & 3.65 & 3.76 & 6337.2 & 0.01 & (10) GOpt & 88.5 & (20) GOpt & 41.7\\
case240\_wecc\_sad & 76494.70& \textbf{5.27} & -- & -- & 10800.0 & 2.61 & (8) \ 2.41 & TO & (16) \ 2.45 & TO \\
case300\_ieee\_sad&16893.92 & 1.09 & 0.13 & 0.10 & 10800.0 & GOpt & (10) GOpt & -- & (20) GOpt & --\\
\midrule
case5\_bgm\_nco & 1082.33 & 10.17 & -- & -- & 6.93 & GOpt & (6) GOpt & -- & (8) GOpt & --  \\
case9\_bgm\_nco & 3087.84 & 10.84 & -- & -- & 2.72 & 10.80   & (16) 10.12 & 7226.6 & (20) 10.12 & 3017.1\\
case9\_na\_cao\_nco & -212.43 & \textbf{-14.97} & -- & -- & 2.14 & -6.90 & (6) -0.04 & 58.6 & (16) -0.73 & 31.4   \\
case9\_nb\_cao\_nco & -247.42 & -15.59 & -- & -- & 2.62 & -6.93 & (6) -0.09 & 31.7 & (8) -0.95 & 133.6 \\
case14\_s\_cao\_nco & 9670.44 & 3.83 & -- & -- & 16.22 & 2.34 & (6) 0.04 & 9626.2 & (10) \ 0.07 & 4092.5 \\
case39\_1\_bgm\_nco & 11221.00 & 3.72 & -- & -- & 103.98 & 3.58 & (8) 3.46 & 14191.8 & (10) \ 3.37 & 17090.3  \\
\bottomrule
\end{tabular}
\end{table*}%

% \input{Sections/NR_arx}
%!TEX root = ../pscc2018.tex

%+++++++++++++++++++++++++++;
%  Section: Conclusions      ;
%+++++++++++++++++++++++++++;
\section{Conclusions}
\label{Sec:conc}
This paper considered the classic ACOPF problem in polar form and developed efficient formulations and algorithms to solve it to global optimality. The key developments in this paper were a) Using state-of-the-art QC relaxations in combination with improved relaxations for trilinear functions based on the convex-hull representation, b) Developing novel mathematical formulations for tight piecewise convex relaxations of trilinear and quadratic functions, and c) Leveraging these tight formulations for global optimization using an adaptive multivariate partitioning approach in combination with bound tightening and effective heuristic branching strategies. Except in a few challenging instances, these methodologies helped to close the optimality gaps of many hard instances. Future directions include testing with other equivalent ACOPF formulations, exploitation of graph sparsitity, and better branching and pruning strategies for piecewise formulations.

%==============================;
%  Include all the references  ;
%==============================;
%\vspace{-0.5cm}
\bibliographystyle{IEEEtran}
\bibliography{references.bib}

\begin{thebibliography}{10}
\providecommand{\url}[1]{#1}
\csname url@samestyle\endcsname
\providecommand{\newblock}{\relax}
\providecommand{\bibinfo}[2]{#2}
\providecommand{\BIBentrySTDinterwordspacing}{\spaceskip=0pt\relax}
\providecommand{\BIBentryALTinterwordstretchfactor}{4}
\providecommand{\BIBentryALTinterwordspacing}{\spaceskip=\fontdimen2\font plus
\BIBentryALTinterwordstretchfactor\fontdimen3\font minus
  \fontdimen4\font\relax}
\providecommand{\BIBforeignlanguage}[2]{{%
\expandafter\ifx\csname l@#1\endcsname\relax
\typeout{** WARNING: IEEEtran.bst: No hyphenation pattern has been}%
\typeout{** loaded for the language `#1'. Using the pattern for}%
\typeout{** the default language instead.}%
\else
\language=\csname l@#1\endcsname
\fi
#2}}
\providecommand{\BIBdecl}{\relax}
\BIBdecl

\bibitem{cain2012history}
M.~B. Cain, R.~P. O’neill, and A.~Castillo, ``History of optimal power flow
  and formulations,'' \emph{Federal Energy Regulatory Commission}, pp. 1--36,
  2012.

\bibitem{carpentier1962contribution}
J.~Carpentier, ``Contribution to the economic dispatch problem,''
  \emph{Bulletin de la Societe Francoise des Electriciens}, vol.~3, no.~8, pp.
  431--447, 1962.

\bibitem{mccormick1976computability}
G.~P. McCormick, ``Computability of global solutions to factorable nonconvex
  programs: Part {I}-- convex underestimating problems,'' \emph{Mathematical
  programming}, vol.~10, no.~1, pp. 147--175, 1976.

\bibitem{hijazi2017convex}
H.~Hijazi, C.~Coffrin, and P.~Van~Hentenryck, ``Convex quadratic relaxations
  for mixed-integer nonlinear programs in power systems,'' \emph{Mathematical
  Programming - C}, vol.~9, no.~3, pp. 321--367, 2017.

\bibitem{lavaei2012zero}
J.~Lavaei and S.~H. Low, ``Zero duality gap in optimal power flow problem,''
  \emph{IEEE Trans. on Power Systems}, vol.~27, no.~1, pp. 92--107, 2012.

\bibitem{kocuk2017matrix}
B.~Kocuk, S.~S. Dey, and X.~A. Sun, ``Matrix minor reformulation and
  {SOCP}-based spatial branch-and-cut method for the {AC} optimal power flow
  problem,'' \emph{arXiv preprint arXiv:1703.03050}, 2017.

\bibitem{coffrin2014nesta}
C.~Coffrin, D.~Gordon, and P.~Scott, ``{NESTA}, the {NICTA} energy system test
  case archive,'' \emph{arXiv preprint arXiv:1411.0359}, 2014.

\bibitem{gopalakrishnan2012global}
A.~Gopalakrishnan, A.~U. Raghunathan, D.~Nikovski, and L.~T. Biegler, ``Global
  optimization of optimal power flow using a branch \& bound algorithm,'' in
  \emph{50th Annual Allerton Conference on Communication, Control, and
  Computing}.\hskip 1em plus 0.5em minus 0.4em\relax IEEE, 2012, pp. 609--616.

\bibitem{nagarajan2016tightening}
H.~Nagarajan, M.~Lu, E.~Yamangil, and R.~Bent, ``Tightening {McC}ormick
  relaxations for nonlinear programs via dynamic multivariate partitioning,''
  in \emph{International Conference on Principles and Practice of Constraint
  Programming}.\hskip 1em plus 0.5em minus 0.4em\relax Springer, 2016, pp.
  369--387.

\bibitem{nagarajan2017adaptive}
H.~Nagarajan, M.~Lu, S.~Wang, R.~Bent, and K.~Sundar, ``An adaptive,
  multivariate partitioning algorithm for global optimization of nonconvex
  programs,'' \emph{arXiv preprint arXiv:1707.02514}, 2017.

\bibitem{coffrin2015strengthening}
C.~Coffrin, H.~L. Hijazi, and P.~Van~Hentenryck, ``Strengthening convex
  relaxations with bound tightening for power network optimization,'' in
  \emph{International Conference on Principles and Practice of Constraint
  Programming}.\hskip 1em plus 0.5em minus 0.4em\relax Springer, 2015, pp.
  39--57.

\bibitem{chen2016bound}
C.~Chen, A.~Atamt{\"u}rk, and S.~S. Oren, ``Bound tightening for the
  alternating current optimal power flow problem,'' \emph{IEEE Trans. on Power
  Systems}, vol.~31, no.~5, pp. 3729--3736, 2016.

\bibitem{wu2017adaptive}
F.~Wu, H.~Nagarajan, A.~Zlotnik, R.~Sioshansi, and A.~M. Rudkevich, ``Adaptive
  convex relaxations for gas pipeline network optimization,'' in \emph{American
  Control Conference, 2017}.\hskip 1em plus 0.5em minus 0.4em\relax IEEE, 2017,
  pp. 4710--4716.

\bibitem{hasan2010piecewise}
M.~Hasan and I.~Karimi, ``Piecewise linear relaxation of bilinear programs
  using bivariate partitioning,'' \emph{AIChE journal}, vol.~56, no.~7, pp.
  1880--1893, 2010.

\bibitem{coffrin2016qc}
C.~Coffrin, H.~L. Hijazi, and P.~Van~Hentenryck, ``The qc relaxation: A
  theoretical and computational study on optimal power flow,'' \emph{IEEE
  Trans. on Power Systems}, vol.~31, no.~4, pp. 3008--3018, 2016.

\bibitem{meyer2004trilinear}
C.~A. Meyer and C.~A. Floudas, ``Trilinear monomials with mixed sign domains:
  Facets of the convex and concave envelopes,'' \emph{Journal of Global
  Optimization}, vol.~29, no.~2, pp. 125--155, 2004.

\bibitem{Rikun1997}
A.~D. Rikun, ``{A Convex Envelope Formula for Multilinear Functions},''
  \emph{Journal of Global Optimization}, vol.~10, pp. 425--437, 1997.

\bibitem{padberg2000approximating}
M.~Padberg, ``Approximating separable nonlinear functions via mixed zero-one
  programs,'' \emph{Operations Research Letters}, vol.~27, no.~1, pp. 1--5,
  2000.

\bibitem{vielma2015mixed}
J.~P. Vielma, ``Mixed integer linear programming formulation techniques,''
  \emph{SIAM Review}, vol.~57, no.~1, pp. 3--57, 2015.

\bibitem{puranik2017domain}
Y.~Puranik and N.~V. Sahinidis, ``Domain reduction techniques for global {NLP}
  and {MINLP} optimization,'' \emph{Constraints}, pp. 1--39, 2017.

\bibitem{bienstock2015strong}
D.~Bienstock and A.~Verma, ``Strong {NP}-hardness of {AC} power flows
  feasibility,'' \emph{arXiv preprint arXiv:1512.07315}, 2015.

\bibitem{lee2001polyhedral}
J.~Lee and D.~Wilson, ``Polyhedral methods for piecewise-linear functions i:
  the lambda method,'' \emph{Discrete applied mathematics}, vol. 108, no.~3,
  pp. 269--285, 2001.

\bibitem{dunning2017jump}
I.~Dunning, J.~Huchette, and M.~Lubin, ``Ju{MP}: A modeling language for
  mathematical optimization,'' \emph{SIAM Review}, vol.~59, no.~2, pp.
  295--320, 2017.

\end{thebibliography}

% Use this to place sponsorships
% \thanksto{Applicable sponsors, if any, should be placed using the \emph{thanksto} command}

\end{document}